\documentclass[10pt]{article}

\usepackage{latexsym,amsfonts}
\usepackage{amssymb,amsthm,upref,amscd}
\usepackage[T1]{fontenc}
\usepackage{times}

\newtheorem{Theorem}{Theorem}[section]
\newtheorem{definition}[Theorem]{Definition}
\newtheorem{lemma}[Theorem]{Lemma}
\newtheorem{proposition}[Theorem]{Proposition}

\newtheorem{Open Problem}[Theorem]{Open Problem}

\makeatletter \@addtoreset{equation}{section} \makeatother
\begin{document}

\title{\bf Existence of a nontrivial solution for a strongly  indefinite   periodic
Schr\"odinger-Poisson  system
 }

\author{Shaowei Chen \thanks{ E-mail address:
chensw@amss.ac.cn (Shaowei Chen)}    \quad \quad Liqin Xiao
  \\ \\
\small \small\it School of Mathematical Sciences, Huaqiao
University, \\
\small Quanzhou  362021,  China\\ }

\date{}
\maketitle
\begin{minipage}{13cm}
{\small {\bf Abstract:}  We consider the Schr\"odinger-Poisson
system
\begin{eqnarray}
\left\{
\begin{array}
[c]{ll}
-\Delta u+V( x) u+|u|^{p-2}u=\lambda \phi u , & \mbox{in }\mathbb{R}^{3},\\
-\Delta\phi= u^{2}, & \mbox{in }\mathbb{R}^{3}.
\end{array}
\right.\nonumber
\end{eqnarray}
where  $\lambda>0$ is a parameter, $3< p<6$, $V\in C(
\mathbb{R}^{3}) $ is $1$-periodic in $x_j$ for $j = 1,2,3$ and
  0 is in a spectral gap of  the operator $-\Delta+V$. This system is strongly indefinite, i.e., the operator $-\Delta+V$
  has infinite-dimensional negative and positive spaces
  and it has  a competitive interplay of the nonlinearities $|u|^{p-2}u$ and $\lambda \phi u$.
  Moreover, the functional corresponding to this system
  does not satisfy the Palai-Smale condition.  Using a new infinite-dimensional linking theorem,
  we prove that, for sufficiently small $\lambda>0,$ this system has a nontrivial solution.  
}
\end{minipage}

\section{Introduction and statement of results}\label{diyizhang}

In this paper, we consider the Schr\"odinger-Poisson system
\begin{eqnarray}\label{nbgftfrr}
\left\{
\begin{array}
[c]{ll}
-\Delta u+V( x) u+|u|^{p-2}u=\lambda \phi u , & \mbox{in }\mathbb{R}^{3},\\
-\Delta\phi= u^{2}, & \mbox{in }\mathbb{R}^{3}.
\end{array}
\right.
\end{eqnarray}
where $3< p<6$ and $\lambda>0$ is a parameter. The potential
function $V$ is a continuous function that is $1$-periodic in
$x_j$ for $j = 1,2,3$. Under this assumption,  $\sigma(L)$, the
spectrum of the operator
\begin{eqnarray}\label{ncxbuc77dy} L=-\Delta+V:
L^{2}(\mathbb{R}^{3})\rightarrow
L^{2}(\mathbb{R}^{3}),\end{eqnarray} is a purely continuous
spectrum that is bounded below and consists of closed disjoint
intervals (\cite[Theorem XIII.100]{Reed2}). Thus, the complement
$\mathbb{R} \setminus \sigma(L)$ consists of  open intervals
called spectral gaps. We assume
\begin{description}
\item{$(\bf{v}).$} $V\in C( \mathbb{R}^{3}) $ is 1-periodic in
$x_j$ for $j = 1,2,3$ and
  0 is in a spectral gap $(-\alpha,\beta)$ of  $-\Delta+V$,
  where $0<\alpha,\beta<+\infty$.
\end{description}

A solution $(u,\phi)$ of (\ref{nbgftfrr}) is called nontrivial if
$u\not\equiv0$ and $\phi\not\equiv 0.$ Our main result is as
follows:

\begin{Theorem}\label{th1}
Suppose that $3<p<6$ and $\bf ( v) $ is satisfied. Then, there
exists $\lambda_0>0$ such that, for any $0<\lambda<\lambda_0,$ the
problem (\ref{nbgftfrr}) has a nontrivial solution.
\end{Theorem}

The Schr\"odinger-Poisson systems such as (\ref{nbgftfrr}) arise
in quantum mechanics, and are related to the study of the
nonlinear Schr\"odinger equation for a particle
 in an electromagnetic field or the Hartree-Fock equation. Such
 systems  have
attracted much attention in recent
 years.
  For
 example, in \cite{lions} and \cite{Lions1}, the system
 \begin{eqnarray}
 -\Delta u+ \omega u+u^{\frac{7}{3}}= \phi u,\
 -\Delta\phi=u^2\quad \mbox{in}\ \mathbb{R}^3\nonumber
 \end{eqnarray}
is used to describe a Hartree model for crystals, where $\omega$
is a positive constant. Problem (\ref{nbgftfrr}) can also be seen
as a nonlinear perturbation of the so-called Choquard system
 \begin{eqnarray}
 -\Delta u+ \omega u=\lambda \phi u,\
 -\Delta\phi=u^2\quad \mbox{in}\ \mathbb{R}^3.\nonumber
 \end{eqnarray}
This system was introduced as an approximation to the Hartree-Fock
theory for a one-component plasma (\cite{ackermann, clapp},
\cite{lieb77}-\cite{Schaftingen}). Problem (\ref{nbgftfrr}) is
also related to
 the so-called
Schr\"odinger-Poisson-Slater system:
$$-\Delta u+V(x)u+\lambda\phi u=f(x,u),\ -\Delta\phi=u^2\quad \mbox{in}\ \mathbb{R}^3$$
where $\lambda>0.$ In recent years, there have been numerous
studies of such systems under various assumptions on $V$ and $f$.
One can, for example, see \cite{alves, Ambrosetti, Azz},
\cite{benci1, tang, Coc, kik,
 Ruiz} for the case $\inf_{\mathbb{R}^3}V
>0$ and \cite{Liu, shaowei, shaowei2} for the case where $V$
changes sign in $\mathbb{R}^3$.  Finally, we should mention that,
in \cite{mugnai}, the author considers the system
$$-\Delta u+\omega u+W'_u(x,u)=\lambda u\phi,\ -\Delta\phi=u^2\quad \mbox{in}\ \mathbb{R}^3,$$
where $\omega$ and $\lambda$ are positive constants. The potential
$W(x,u)$ is  nonnegative,  radially symmetric in $x$ and, even in
$u$ and satisfies some growth conditions. The model nonlinearity
is $\frac{1}{p}|u|^p$ with $2<p<6.$ Using the symmetric mountain
pass theorem (see \cite{AR}), the author obtained infinitely many
radially symmetric solutions of this system.

 From a mathematical point of view, problem
(\ref{nbgftfrr})
 possesses some interesting properties. First of all, this problem
has a variational structure. Let
$\mathcal{D}^{1,2}(\mathbb{R}^{3})$ be the Hilbert space
$$\mathcal{D}^{1,2}(\mathbb{R}^{3})=\{u\in L^6(\mathbb{R}^3)\
 |\ \int_{\mathbb{R}^3}|\nabla u|^2dx<\infty\}$$ with inner
 product $$(u,v)_{\mathcal{D}^{1,2}(\mathbb{R}^{3})}=\int_{\mathbb{R}^3}\nabla u\nabla vdx.$$
For $(u,\phi)\in H^1(\mathbb{R}^3) \times
\mathcal{D}^{1,2}(\mathbb{R}^3)$, let
$$J_{\lambda}(u,\phi)=\frac{1}{2}\int_{\mathbb{R}^3}(|\nabla u|^2+V(x)u^2)dx
+\frac{\lambda}{4}\int_{\mathbb{R}^3}|\nabla\phi|^2dx+\frac{1}{p}\int_{\mathbb{R}^3}|u|^pdx
-\frac{\lambda}{2}\int_{\mathbb{R}^3}u^2\phi dx.$$ Then, the
critical points of $J_{\lambda}$ solve (\ref{nbgftfrr}). Because
$V$ is $1$-periodic in every variable, this functional is
invariant under a $1$-periodic translation. As a consequence, the
functional does not satisfy   the Palais-Smale condition (see
\cite{Willem}). Second, because $0$ is in a spectral gap of
$-\Delta+V$, the quadratic form $\int_{\mathbb{R}^3}(|\nabla
u|^2+V(x)u^2)dx$ has infinite-dimensional negative  and positive
spaces. This case is called strongly indefinite. Finally, since
$J_\lambda$ is strongly indefinite, it is natural to use the
infinite-dimensional linking theorem (see \cite[Theorem 3.4]{KS}
or \cite[Theorem 6.10]{Willem})  to obtain critical points of
$J_\lambda$. However, because of  the competitive interplay of the
nonlinearities $\frac{1}{p}\int_{\mathbb{R}^3}|u|^pdx$ and
$\lambda\int_{\mathbb{R}^3}u^2\phi dx$, $J_\lambda$ does not
satisfy   the global linking condition (see (6.4) in
\cite{Willem}) or the $\tau$-upper semi-continuous assumption (see
(6.3) in \cite{Willem}) in this theorem. To the best of our
knowledge,  variational problems possessing all these properties
have never been studied before. Our study of this problem will
shed some light on other variational problems possessing
properties similar to those mentioned above.

In this paper,  we first,   modify the nonlinear terms of problem
(\ref{nbgftfrr}) such that the variational functional for the
modified system of equations has a global infinite-dimensional
linking structure. We then use a new infinite-dimensional linking
theorem (see Theorem \ref{b99d443} in the appendix) to obtain a
nontrivial solution for the modified system. This theorem replaces
the $\tau$-upper semi-continuous assumption in the classical
infinite-dimensional linking theorem (see \cite[Theorem
6.10]{Willem} or \cite[Theorem 3.4]{KS}) with other assumption
(see (\ref{nx99s8s7yy}) in Theorem \ref{b99d443}). For  reader's
convenience, we give the proof of this theorem in the appendix.
Finally, we use a blow-up argument to show that if $\lambda>0$ is
sufficiently small, then the nontrivial solution of the modified
system is in fact a nontrivial solution of (\ref{nbgftfrr}).

 \medskip

 \noindent{\bf Notation.} $B_r(a)$ denotes the  open ball of radius $r$ and center $a$.
For a Banach space $X,$ we denote the dual space of $X$ by $X'$,
and denote   strong and   weak convergence in $X$  by
$\rightarrow$ and $\rightharpoonup$, respectively. For $\varphi\in
C^1(X;\mathbb{R}),$ we denote the Fr\'echet derivative of
$\varphi$ at $u$ by $\varphi'(u)$.  The Gateaux derivative of
$\varphi$ is denoted by $\langle \varphi'(u), v\rangle,$ $\forall
u,v\in X.$ $L^q(\mathbb{R}^3)$ and $L^q_{loc}(\mathbb{R}^3)$
denote the standard $L^q$ space and the locally $q$-integrable
function space, respectively $(1\leq q\leq\infty)$, and
$H^1(\mathbb{R}^3)$ denotes the standard Sobolev space with norm
$||u||_{H^1}=(\int_{\mathbb{R}^3}(|\nabla u|^2+u^2)dx)^{1/2}.$ Let
$\Omega$ be a domain in $\mathbb{R}^N$ ($N\geq 1$).
$C^\infty_0(\Omega)$ is the space of infinitely differentiable
functions with compact support in $\Omega$.

\section{A modified  system for (\ref{nbgftfrr})}

Let $\eta\in C^\infty_0(\mathbb{R})$ be such that
$$0\leq\eta\leq 1;\ \eta(t)=1,\ r\in [-1,1]; \ \eta(t)=0,\ |t|\geq 2;
\ |\eta'|\leq 2;\ \eta(-t)=\eta(t),\ \forall t\in\mathbb{R}. $$
 For $n\in\mathbb{N}$ and $t\in\mathbb{R}$, let
  $$\eta_n(t)=\eta(t/n),\ F_n(t)=t^2\eta_n(t)\ \mbox{and}\ f_n(t)=F'_n(t)/2.$$

 Consider the following system
\begin{equation}\label{nbgftfra44sr}
\left\{
\begin{array}
[c]{ll}
-\Delta u+V( x) u+|u|^{p-2}u=\lambda\phi f_n(u) , & \mbox{in }\mathbb{R}^{3},\\
-\Delta\phi=F_n(u), & \mbox{in }\mathbb{R}^{3}.
\end{array}
\right.
\end{equation}

We can prove that the second equation has a unique solution
$\phi_u$. Substituting $\phi_u$ into the first equation of
(\ref{nbgftfra44sr}),  the problem can be transformed into a one
variable equation. In fact, we  have the following lemma:

\begin{lemma}\label{fhyr66rtfre}
For any $u \in H^1(\mathbb{R}^3)$, there exists a unique
$\phi_{n,u} \in \mathcal{D}^{1,2}(\mathbb{R}^3)$ that is a
solution of
\begin{eqnarray}\label{nx88fdhyftt}
-\Delta\phi=F_n(u),\quad \mbox{in}\ \mathbb{R}^3. \end{eqnarray}
Moreover,
\begin{itemize}
\item[{(i).}]
$\phi_{n,u}(x)=\int_{\mathbb{R}^3}\frac{F_n(u(y))}{4\pi|x-y|}dy,$
$x\in\mathbb{R}^3.$ \item[{(ii).}]
 $\int_{\mathbb{R}^3}|\nabla\phi_{n,u}|^2dx=\int_{\mathbb{R}^3}F_n(u)\phi_{n,u}dx
 =\int_{\mathbb{R}^3}\int_{\mathbb{R}^3}\frac{F_n(u(x))F_n(u(y))}{4\pi|x-y|}dxdy$.
\item[{(iii).}] Let $$C_*=\inf_{0\neq v\in
\mathcal{D}^{1,2}(\mathbb{R}^3)}\frac{\int_{\mathbb{R}^3}|\nabla
v|^2dx}{(\int_{\mathbb{R}^3}|v|^6dx)^{1/3}}>0$$ be the Sobolev
constant (see, for example, \cite[Theorem 1.8]{Willem}). Then
$$\int_{\mathbb{R}^3}|\nabla\phi_{n,u}|^2dx\leq C^{-1}_*(\int_{\mathbb{R}^3}|u|^{\frac{12}{5}}dx)^{\frac{5}{3}},\ \forall u\in H^1(\mathbb{R}^3).$$
\item[{(iv).}]There exists a positive constant $D$, which is
independent of $n$ and $u$, such that $0\leq\phi_{n,u}\leq Dn^2$
in $\mathbb{R}^3$.
\end{itemize}
\end{lemma}
\noindent{\bf Proof.} Because $0\leq F_n(t)\leq t^2$ for all $n$
and $t$,  the H\"older and the Sobolev inequalities implies that,
for any $u\in H^1(\mathbb{R}^3)$ and $v\in
\mathcal{D}^{1,2}(\mathbb{R}^3)$,
\begin{eqnarray}\label{ndh878fryf}
|\int_{\mathbb{R}^3}F_n(u)vdx|
&\leq&(\int_{\mathbb{R}^3}|F_n(u)|^{\frac{6}{5}}dx)^{\frac{5}{6}}(\int_{\mathbb{R}^3}|v|^6dx)^{\frac{1}{6}}\nonumber\\
&\leq&(\int_{\mathbb{R}^3}|u|^{\frac{12}{5}}dx)^{\frac{5}{6}}(\int_{\mathbb{R}^3}|v|^6dx)^{\frac{1}{6}}
\leq C^{-\frac{1}{2}}_*
(\int_{\mathbb{R}^3}|u|^{\frac{12}{5}}dx)^{\frac{5}{6}}(\int_{\mathbb{R}^3}|\nabla
v|^2dx)^{\frac{1}{2}}.
\end{eqnarray}
 It follows that, for fixed $u\in H^1(\mathbb{R}^3)$, $v\in
\mathcal{D}^{1,2}(\mathbb{R}^3)\mapsto
\int_{\mathbb{R}^3}F_n(u)vdx$ is a bounded linear functional in
$\mathcal{D}^{1,2}(\mathbb{R}^3).$ Then, by the Riesz theorem,
there exists a unique $\phi_{n,u}\in
\mathcal{D}^{1,2}(\mathbb{R}^3)$ such that
\begin{eqnarray}\label{nchyfteee}
\int_{\mathbb{R}^3}\nabla \phi_{n,u}\nabla
vdx=\int_{\mathbb{R}^3}F_n(u)vdx,\ \forall v\in
\mathcal{D}^{1,2}(\mathbb{R}^3).
 \end{eqnarray} It follows that
$\phi_{n,u}$ is the unique solution of (\ref{nx88fdhyftt}). Then,
by the theory of Poisson's equation (see, for example,  Theorem
2.2.1 of \cite{Evans}), we obtain the expression for $\phi_{n,u}$
in $(i)$. Moreover, choosing $v=\phi_{n,u}$ in (\ref{nchyfteee}),
we obtain  result $(ii)$ of this lemma.

Choosing $v=\phi_{n,u}$ in (\ref{ndh878fryf}) and using the first
equality in  $(ii)$, we obtain  result $(iii).$

Since $\phi_{n,u}$ is the  solution of (\ref{fhyr66rtfre}),  the
regularity theory for elliptic equations (see, for example,
\cite[Theorem 8.17]{GT}) implies that  there exists a positive
constant $C_1$ that is independent of
 $n$, $u$ and $y\in\mathbb{R}^3$, such that, for any
$y\in\mathbb{R}^3$,
\begin{eqnarray}\label{nnbvyyfh765}
||\phi_{n,u}||_{L^\infty(B_1(y))}\leq
C_1(\int_{B_2(y)}|F_n(u)|^2)^{\frac{1}{2}}.
\end{eqnarray}
 From the definition  of $F_n$,  we have
$0\leq F_n(t)\leq (2 n)^2,$ $\forall t\in\mathbb{R}.$ It follows
that
\begin{eqnarray}\label{mvcn88fuyf}
C_1(\int_{B_2(y)}|F_n(u)|^2)^{\frac{1}{2}}\leq
C_1(\int_{B_2(y)}(2n)^4
)^{\frac{1}{2}}=4C_1(\int_{B_2(0)}dx)^{1/2}n^2.
\end{eqnarray}
Choosing $D=4C_1(\int_{B_2(0)}dx )^{\frac{1}{2}},$ by
(\ref{nnbvyyfh765}) and (\ref{mvcn88fuyf}), we get that
$\phi_{n,u}\leq Dn^2 $ in $\mathbb{R}^3$. Finally, since
$F_n(u)\geq 0$ in $\mathbb{R}^3$,  the maximum principle (see, for
example, \cite{GT}) implies that $\phi_{n,u}\geq 0$ in
$\mathbb{R}^3$. \hfill$\Box$

\bigskip

For $u\in H^1(\mathbb{R}^3)$, let
\begin{eqnarray}\label{ncxbttdref}
&&\Phi_{n,\lambda}(u)\nonumber\\
&=&\frac{1}{2}\int_{\mathbb{R}^3}(|\nabla u|^2+V(x)u^2)dx
+\frac{1}{p}\int_{\mathbb{R}^3}|u|^pdx
-\frac{\lambda}{4}\int_{\mathbb{R}^3}F_n(u)\phi_{n,u}
dx\nonumber\\
&=&\frac{1}{2}\int_{\mathbb{R}^3}(|\nabla u|^2+V(x)u^2)dx
+\frac{1}{p}\int_{\mathbb{R}^3}|u|^pdx
-\frac{\lambda}{4}\int_{\mathbb{R}^3}\int_{\mathbb{R}^3}\frac{F_n(u(x))F_n(u(y))}{4\pi|x-y|}dxdy.
\end{eqnarray}
A direct computation shows that the derivative of
$\Phi_{n,\lambda}$ is
\begin{eqnarray}\label{bc88dytr22}
\langle \Phi'_{n,\lambda}(u), v\rangle
&=&\int_{\mathbb{R}^3}(\nabla u\nabla v+V(x)uv)dx
+\int_{\mathbb{R}^3}|u|^{p-2}uvdx\nonumber\\
&&
-\lambda\int_{\mathbb{R}^3}\int_{\mathbb{R}^3}\frac{f_n(u(x))v(x)F_n(u(y))}{4\pi|x-y|}dxdy,\
\forall u,v\in H^1(\mathbb{R}^3).
\end{eqnarray}

From   $(i)$ of Lemma \ref{fhyr66rtfre},  and (\ref{bc88dytr22}),
we have the following

\begin{lemma}\label{vd88dytrdf11}
The following statements are equivalent:
\begin{itemize}
\item[{(i).}]  $(u, \phi) \in
H^1(\mathbb{R}^3)\times\mathcal{D}^{1,2}(\mathbb{R}^3)$ is a
solution  of (\ref{nbgftfra44sr}).
 \item[{(ii).}]  $u$ is a
critical point of $\Phi_{n,\lambda}$ and $\phi=\phi_{n,u}$.
\end{itemize}
\end{lemma}

\section{Existence of a nontrivial  solution for  (\ref{nbgftfra44sr})}
Recall that $L$ is the operator defined by (\ref{ncxbuc77dy}).  We
denote by $|L|^{1/2}$ the square root of the absolute value of
$L.$ The domain of $|L|^{1/2} $ is the space
$$X:=H^{1}(\mathbb{R}^{3}).$$ On $X,$ we choose the inner product
$(u,v)=\int_{\mathbb{R}^{3}}|L|^{1/2}u\cdot |L|^{1/2}v dx$ and the
corresponding norm $||u||=\sqrt{ (u,u)}.$ Since  $0$ lies in a gap
of the essential spectrum of $L,$ there exists an orthogonal
decomposition $X=Y\oplus Z$ such that $Z$ and $Y$ are the positive
and negative spaces corresponding to the spectral decomposition of
$L$. Since $V$ is $1$-periodic for all variables, they are
invariant under the action of $\mathbb{Z}^{3},$ i.e., for any
$u\in Y$ or $u\in Z$ and for any $\mathbf{k}=(n_{1},n_2,n_{3})\in
\mathbb{Z}^{3},$ $u(\cdot-\mathbf{k})$ is also in $Y$ or $Z.$
Furthermore,
\begin{equation}\label{hnbxxc5r8871z} \forall u\in Y,\
\int_{\mathbb{R}^{3}}(|\nabla
u|^{2}+Vu^{2})dx=(u,u)=||u||^{2},\end{equation}
\begin{equation}\label{un77715490kjh76}
\forall u\in Z,\ \int_{\mathbb{R}^{3}}(|\nabla
u|^{2}+Vu^{2})dx=-(u,u)=-||u||^{2}.
\end{equation}
 Let $Q:X\rightarrow Z$, $P:X\rightarrow Y$ be the
orthogonal projections. By (\ref{hnbxxc5r8871z}) and
(\ref{un77715490kjh76}),
\begin{eqnarray}\label{vd66dtrdpp}
\int_{\mathbb{R}^{3}}(|\nabla u|^{2}+Vu^{2})dx=||Pu||^2-||Qu||^2,\
\forall u\in X.
\end{eqnarray}
From $X=Y\oplus Z$, we get that, for any $u\in X,$
\begin{eqnarray}\label{nfhyfyrttr11}
u=Pu+Qu,\quad ||u||^2=||Pu||^2+||Qu||^2.
\end{eqnarray}

This is the standard  variational setting for the quadratic form
$\int_{\mathbb{R}^3}(|\nabla u|^2+V(x)u^2)dx$. See section 6.4 of
\cite{Willem} for more details.

 By (\ref{vd66dtrdpp}) and (\ref{ncxbttdref}), we have
\begin{eqnarray}\label{ncxbtsdsesqf}
-\Phi_{n,\lambda}(u)=\frac{1}{2}||Qu||^2-\frac{1}{2}||Pu||^2
-\frac{1}{p}\int_{\mathbb{R}^3}|u|^pdx
+\frac{\lambda}{4}\int_{\mathbb{R}^3}F_n(u)\phi_{n,u} dx.
\end{eqnarray}
Moreover, by (\ref{bc88dytr22}), for any $u,v\in X,$
\begin{eqnarray}\label{hdgtdrddd}
&&\langle-\Phi'_{n,\lambda}(u),v\rangle\nonumber\\
&=&(Qu, v)-(Pu,v)-\int_{\mathbb{R}^3}|u|^{p-2}uvdx
+\lambda\int_{\mathbb{R}^3}\int_{\mathbb{R}^3}\frac{f_n(u(x))v(x)F_n(u(y))}{4\pi|x-y|}dxdy\nonumber\\
&=&(Qu, v)-(Pu,v)-\int_{\mathbb{R}^3}|u|^{p-2}uvdx
+\lambda\int_{\mathbb{R}^3}\phi_{n,u}f_n(u) vdx.
\end{eqnarray}

We will prove that if $\lambda>0$ is sufficiently small, then
$-\Phi_{n,\lambda}$  satisfies the global linking condition (see
Lemma \ref{b99d443iu}). However, because of the nonlinearity,
$\frac{\lambda}{4}\int_{\mathbb{R}^3}F_n(u)\phi_{n,u} dx$,
$-\Phi_{n,\lambda}$ does not satisfies the $\tau$-upper continuous
assumption (see (6.3) in \cite{Willem}). Therefore, to obtain
critical points of $-\Phi_{n,\lambda}$, we have to use the new
infinite-dimensional linking theorem (Theorem \ref{b99d443} in
appendix).

 Let $\{e_k\}$ be a
total orthonormal sequence in $Y$ and
\begin{eqnarray}\label{nvcooiugyu}
|||u|||=\max\Big\{||Qu||,\sum^{\infty}_{k=1}\frac{1}{2^{k+1}}|(Pu,e_k)|\Big\}.
\end{eqnarray} For $R > r > 0$ and  $u_0 \in Z$ with $||u_0||
 = 1$,
 set
\begin{eqnarray}
N = \{u \in Z \ |\  ||u||
 = r\},\  M = \{u+tu_0\ |\  u\in Y ,\ t\geq 0,\
||u+tu_0|| \leq R\}\nonumber \end{eqnarray} and
\begin{eqnarray}
\partial M=\{u\in Y\ |\ ||u||\leq
R\}\cup\{u+tu_0\ |\ u\in Y,\ t>0,\ ||u+tu_0||=R \}.\nonumber
\end{eqnarray}

\begin{lemma} \label{b99d443iu}
 The functional $-\Phi_{n,\lambda}$ satisfies the following
\begin{description}  \item {$(a)$}  $-\Phi'_{n,\lambda}$  is weakly
sequentially continuous, where  the weakly sequential continuity
is defined in Theorem \ref{b99d443} in the appendix. \item {$(b)$}
there exist $\delta>0,$ $R
> r
> 0$, $u_0 \in Z  $ with $||u_0||=1$ and $\lambda'_n>0$ for any $n\in\mathbb{N}$
such that if $0<\lambda<\lambda'_n,$ then
\begin{eqnarray}\label{qanx99s8s7yy}
\inf_N (-\Phi_{n,\lambda})>\max\Big\{0,\ \sup_{\partial M}
(-\Phi_{n,\lambda}),
\sup_{|||u|||\leq\delta}(-\Phi_{n,\lambda}(u))\Big\}
\end{eqnarray}
and
\begin{eqnarray}\label{qam99iuxyyxy}
\sup_M(-\Phi_{n,\lambda})<+\infty.
\end{eqnarray}
\end{description}
\end{lemma}
\noindent{\bf Proof.} $(a)$  Let $u\in X$ and $\{u_k\}\subset X$
be such that $u_k\rightharpoonup u$ as $k\rightarrow\infty$.  It
follows that
\begin{eqnarray}\label{hd66dtrffq}
(Qu_k,v)\rightarrow (Qu,v),\ (Pu_k,v)\rightarrow (Pu,v),\
k\rightarrow\infty, \forall v\in X,
\end{eqnarray}
and $u_k\rightarrow u$ in $L^s_{loc}(\mathbb{R}^3)$ for any $1\leq
s<6$.  As consequences, for any $v\in C^\infty_0(\mathbb{R}^3),$
as $k\rightarrow\infty,$
\begin{eqnarray}\label{nfdhhhhh8y}
\int_{\mathbb{R}^3}|u_k|^{p-2}u_kvdx\rightarrow\int_{\mathbb{R}^3}|u|^{p-2}uvdx,
\end{eqnarray}
and
\begin{eqnarray}\label{nc90viufyf}
\int_{\mathbb{R}^3}F_n(u_k)vdx\rightarrow\int_{\mathbb{R}^3}F_n(u)vdx.
\end{eqnarray}

From (\ref{nchyfteee}) and  (\ref{nc90viufyf}), we have that,  for
any $v\in C^\infty_0(\mathbb{R}^3),$ as $k\rightarrow\infty,$
\begin{eqnarray}
\int_{\mathbb{R}^3}\nabla\phi_{n,u_k}\nabla vdx\rightarrow
\int_{\mathbb{R}^3}\nabla\phi_{n,u}\nabla vdx.\nonumber
\end{eqnarray}
This implies $\phi_{n,u_k}\rightharpoonup \phi_{n,u}$ in
$\mathcal{D}^{1,2}(\mathbb{R}^3)$. Consequently,
$\phi_{n,u_k}\rightarrow \phi_{n,u}$ in $L^s_{loc}(\mathbb{R}^3)$
for any $1\leq s<6$. Together with $u_k\rightarrow u$ in
$L^s_{loc}(\mathbb{R}^3)$ for any $1\leq s<6$, this yields that,
for any $v\in C^\infty_0(\mathbb{R}^3),$ as $k\rightarrow\infty,$
\begin{eqnarray}\label{mvn88fyrttr}
\int_{\mathbb{R}^3}\phi_{n,u_k}f_n(u_k)vdx\rightarrow\int_{\mathbb{R}^3}\phi_{n,u}f_n(u)vdx.
\end{eqnarray}
From (\ref{hdgtdrddd}), (\ref{hd66dtrffq}), (\ref{nfdhhhhh8y}) and
(\ref{mvn88fyrttr}), we get that, as $k\rightarrow\infty,$
$$\langle-\Phi'_{n,\lambda}(u_k), v\rangle\rightarrow \langle-\Phi'_{n,\lambda}(u),
v\rangle,\ \forall v\in C^\infty_0(\mathbb{R}^3).$$ Therefore,
$-\Phi'_{n,\lambda}$ is weakly sequentially continuous. Moreover,
$-\Phi_{n,\lambda}$ maps bounded sets into bounded sets, hence
$\sup_{M}(-\Phi_{n,\lambda})<+\infty$.

\medskip

$(b)$ If $u\in Z$, then $Pu=0$ and $Qu=u.$ As $\phi_{n,u}\geq 0$
(see Lemma \ref{fhyr66rtfre}) and   $F_n\geq 0,$ we have
$\int_{\mathbb{R}^3}F_n(u)\phi_{n,u}dx\geq 0$ for all $u\in X$.
Then, using the Sobolev inequality
\begin{eqnarray}\label{m000odiiii1}||u||_{L^p(\mathbb{R}^3)}\leq C'||u||,\end{eqnarray} we get for the
definition of $-\Phi_{n,\lambda}$ that, for any $u\in Z,$
\begin{eqnarray}\label{ncuudyfh}
-\Phi_{n,\lambda}(u)\geq
\frac{1}{2}||u||^2-\frac{1}{p}\int_{\mathbb{R}^3}|u|^pdx\geq
\frac{1}{2}||u||^2-\frac{C'^p}{p}||u||^p.\nonumber
\end{eqnarray}
 Let $r=C'^{-p/(p-2)}$. Then,  for $N=\{u\in Z\ |\ ||u||=r\}$,
\begin{eqnarray}\label{cb99dudy}
\inf_{N}(-\Phi_{n,\lambda})\geq(\frac{1}{2}-\frac{1}{p})C'^{-2p/(p-2)}>0.
\end{eqnarray}

Let $C>0$ be such that
\begin{eqnarray}\label{mciiiicjpp09}
||u||_{L^2(\mathbb{R}^3)}\leq C||u||,\ \forall u\in X.
\end{eqnarray} Let
$$\lambda'_n=(C^2Dn^2)^{-1},$$ where $D$ is the constant appearing
in  Lemma \ref{fhyr66rtfre}$(iv)$.   Then, for any
$0<\lambda\leq\lambda'_n$ and $u\in X,$  $F_n(t)\leq t^2$ and
$\phi_{n,u}\leq Dn^2$ (see  Lemma \ref{fhyr66rtfre}$(iv)$) yield
\begin{eqnarray}\label{bcgttf88uy}
-\Phi_{n,\lambda}(u) &=&
\frac{1}{2}||Qu||^2-\frac{1}{2}||Pu||^2-\frac{1}{p}\int_{\mathbb{R}^3}|u|^pdx
+\frac{\lambda}{4}\int_{\mathbb{R}^3}F_n(u)\phi_{n,u}dx\nonumber\\
&\leq&\frac{1}{2}||Qu||^2-\frac{1}{2}||Pu||^2-\frac{1}{p}\int_{\mathbb{R}^3}|u|^pdx+
\frac{\lambda}{4}\int_{\mathbb{R}^3}u^2\cdot(Dn^2)dx\nonumber\\
&\leq&\frac{1}{2}||Qu||^2-\frac{1}{2}||Pu||^2-\frac{1}{p}\int_{\mathbb{R}^3}|u|^pdx+\frac{1}{4}||u||^2\nonumber\\
&=&\frac{3}{4}||Qu||^2-\frac{1}{4}||Pu||^2-\frac{1}{p}\int_{\mathbb{R}^3}|u|^pdx.
\end{eqnarray}

Let $u_0\in Z$ be such that $||u_0||=1.$  And let $u=v+tu_0\in
Y\oplus \mathbb{R}u_0.$ By (\ref{bcgttf88uy}), we have
\begin{eqnarray}\label{mnghgyuyu00}-\Phi_{n,\lambda}(u)\leq
\frac{3}{4}t^2-\frac{1}{4}||v||^2-\frac{1}{p}\int_{\mathbb{R}^3}|v+tu_0|^pdx.
\end{eqnarray}
From (\ref{mnghgyuyu00}) and the proof of Lemma 6.14 in
\cite{Willem}, we obtain
\begin{eqnarray}\label{vdgttdg112}
-\Phi_{n,\lambda}(u)\rightarrow -\infty,\ \mbox{as}\ \
||u||\rightarrow\infty,\  u\in Y\oplus \mathbb{R}u_0.
\end{eqnarray}
Moreover, for $0<\lambda\leq\lambda'_n$ and $u\in Y,$
(\ref{bcgttf88uy}) implies
\begin{eqnarray}\label{md090d9d8ww}
-\Phi_{n,\lambda}(u)\leq-\frac{1}{4}||u||^2-\frac{1}{p}\int_{\mathbb{R}^3}|u|^pdx\leq
0.\nonumber
\end{eqnarray}
 Together with (\ref{vdgttdg112}), this yields
\begin{eqnarray}\label{mv99fugfytt}
\sup_{\partial M}(-\Phi_{n,\lambda})\leq
0<\inf_N(-\Phi_{n,\lambda}).
\end{eqnarray}

From (\ref{bcgttf88uy}) and the definition of $|||\cdot|||$ (see
(\ref{nvcooiugyu})), we have, for $0<\lambda\leq\lambda'_n,$
\begin{eqnarray}\label{ncx99cifyuhh}
-\Phi_{n,\lambda}(u) \leq
\frac{3}{4}||Qu||^2\leq\frac{3}{4}|||u|||^2.
\end{eqnarray}
Choosing
$\delta=(\frac{2}{3}(\frac{1}{2}-\frac{1}{p})C'^{-2p/(p-2)})^{1/2}$,
 (\ref{ncx99cifyuhh}) and (\ref{cb99dudy}) give that
$$\sup_{|||u|||\leq\delta}(-\Phi_{n,\lambda}(u))<\inf_{N}(-\Phi_{n,\lambda}).$$
Together with (\ref{mv99fugfytt}), this yields
(\ref{qanx99s8s7yy}).\hfill$\Box$

\begin{lemma}\label{hf99ofiuufyf} Let $V_-=\max\{-V,0\}$ and
$M=||V_-||^{\frac{1}{p-2}}_{L^\infty(\mathbb{R}^3)}+2$. For any
$n\in\mathbb{N}$, there exists $\lambda''_n>0$ such that if
$0<\lambda\leq\lambda''_n$ and $\{u_k\}$ is  a $(\overline{C})_c$
sequence
 for $-\Phi_{n,\lambda}$, i.e.
 \begin{eqnarray}\label{m99kllll112}
\sup_n( -\Phi_{n,\lambda}(u_k))\leq c,\
(1+||u_k||)||-\Phi'_{n,\lambda}(u_k)||_{X'}\rightarrow 0,\
\mbox{as}\ k\rightarrow\infty,
\end{eqnarray}
 then
\begin{eqnarray}\label{l009dttdgf}
\int_{\varpi_k}|u_k|^pdx\rightarrow 0,\ k\rightarrow\infty,
\end{eqnarray}
where
$$\varpi_k=\{x\in\mathbb{R}^3\ |\ |u_k(x)|\geq M \}.$$

\end{lemma}
\noindent{\bf Proof.}  Let
$$v_k(x)=\max\{u_k-M+1,0\}.$$
It is easy to see that $\int_{\mathbb{R}^3}|\nabla
v_k|^2dx\leq\int_{\mathbb{R}^3}|\nabla u_k|^2dx$ and
$\int_{\mathbb{R}^3}| v_k|^2dx\leq\int_{\mathbb{R}^3}|u_k|^2dx$.
Therefore, there exists $C''>0$ such that $||v_k||\leq
C''||u_k||$, $\forall k\in\mathbb{N}.$  Then, by
(\ref{m99kllll112}), we get
$\langle\Phi'_{n,\lambda}(u_k),v_k\rangle=o(1),$ here $o(1)$
denotes the infinitesimal depending only on $k,$ i.e.,
$o(1)\rightarrow 0$ as $k\rightarrow \infty.$ Together with
(\ref{bc88dytr22}),  this yields
\begin{eqnarray}\label{bv88dgffd11}
o(1)&=&\int_{\mathbb{R}^3}\nabla u_k\nabla
v_kdx+\int_{\mathbb{R}^3}V(x)u_kv_kdx+\int_{\mathbb{R}^3}|u_k|^{p-2}u_kv_kdx-\lambda\int_{\mathbb{R}^3}f_n(u_k)\phi_{n,u_k}v_kdx\nonumber\\
&=&\int_{\widetilde{\varpi}^+_k}|\nabla
v_k|^2dx+\int_{\widetilde{\varpi}^+_k}V_+(x)
u_kv_kdx\nonumber\\
&&+\int_{\widetilde{\varpi}^+_k}(|u_k|^{p-2}-V_--\lambda u^{-1}_k
f_n(u_k)\phi_{n,u_k})u_kv_kdx
\end{eqnarray}
where $V_+=V+V_-$ and $$\widetilde{\varpi}^+_k=\{x\in\mathbb{R}^3\
|\ u_k(x)\geq M-1 \}.$$ From $0\leq\phi_{n,u_k}\leq Dn^2$ and
$|t^{-1}f_n(t)|\leq 3$, we deduce that if  $$0<\lambda\leq
\frac{1}{6Dn^2}:=\lambda''_n,$$ then, for any
$x\in\widetilde{\varpi}_k,$
$$|u_k|^{p-2}-V_--\lambda u^{-1}_k
f_n(u_k)\phi_{n,u_k}> 0.$$ Together with (\ref{bv88dgffd11}) and
the fact that $V_+,$ $u_k$ and $v_k$ are nonnegative on
$\widetilde{\varpi}_k$, this implies
$$\int_{\widetilde{\varpi}^+_k}|\nabla
v_k|^2dx=o(1).$$ Then, by the Sobolev inequality, we obtain
\begin{eqnarray}
\int_{\widetilde{\varpi}^+_k}| v_k|^6dx=o(1).\nonumber
\end{eqnarray}
Let$$\varpi^+_k=\{x\in\mathbb{R}^3\ |\ u_k(x)\geq M \}.$$ Then
$\varpi^+_k\subset \widetilde{\varpi}^+_k$. And on $\varpi^+_k$,
$v^6_k\geq u^6_k/M^6\geq u^p_k/M^p$. It follows that
\begin{eqnarray}\label{ncggdfdfdkm}
\int_{\varpi^+_k}| u_k|^pdx\leq M^p\int_{\widetilde{\varpi}^+_k}|
v_k|^6dx=o(1).
\end{eqnarray}
Similarly, we have
\begin{eqnarray}\label{1ncggdfdfdkm}
\int_{\varpi^-_k}| u_k|^pdx=o(1).
\end{eqnarray}
where $\varpi^-_k=\{x\in\mathbb{R}^3\ |\ -u_k(x)\geq M \}.$
(\ref{l009dttdgf}) follows from (\ref{ncggdfdfdkm}) and
(\ref{1ncggdfdfdkm}) immediately.\hfill$\Box$

\begin{lemma}\label{fbf99fifio}
For any $n\in\mathbb{N}$, there exists $\lambda'''_n>0$ such that,
if $0<\lambda<\lambda'''_n$ and $\{u_k\}$ is  a $(\overline{C})_c$
sequence
 for $-\Phi_{n,\lambda}$,  then $\{u_k\}$ is bounded in $X.$
\end{lemma}
\noindent{\bf Proof.} From
$(1+||u_k||)||-\Phi'_{n,\lambda}(u_k)||_{X'}\rightarrow 0$, we
have  $$\langle -\Phi'_{n,\lambda}(u_k), Qu_k\rangle=o(1)\quad
\mbox{and}\quad \langle -\Phi'_{n,\lambda}(u_k),
Pu_k\rangle=o(1),$$ where $o(1)$ denotes the infinitesimal
depending only on $k,$ i.e., $o(1)\rightarrow 0$ as $k\rightarrow
\infty.$ Together with (\ref{hdgtdrddd}), this yields
\begin{eqnarray}\label{hf999f876r6rr}
||Qu_k||^2=\int_{\mathbb{R}^3}|u_k|^{p-2}u_k\cdot
Qu_kdx-\lambda\int_{\mathbb{R}^3}f_n(u_k)\phi_{n,u_k}\cdot
Qu_kdx+o(1)
\end{eqnarray}
and
\begin{eqnarray}\label{hfuu7g7tyyt1}
||Pu_k||^2=-\int_{\mathbb{R}^3}|u_k|^{p-2}u_k\cdot
Pu_kdx+\lambda\int_{\mathbb{R}^3}f_n(u_k)\phi_{n,u_k}\cdot
Pu_kdx+o(1).
\end{eqnarray}
  For $\epsilon>0,$ let $$A_{\epsilon,k}=\{x\in\mathbb{R}^3\ |\ |u_k(x)|<\epsilon\}.$$
  And recall that $\varpi_k$ is the set defined in Lemma
 \ref{hf99ofiuufyf}.
 Using  $0\leq\phi_{n,u_k}\leq Dn^2$, $| f_n(t)|\leq 5t$ and $|u_k|\leq M$ on $\mathbb{R}^3\setminus\varpi_k$,
  we get from
 (\ref{hf999f876r6rr}) that
\begin{eqnarray}\label{ijnwdvb66}
&&||Qu_k||^2\nonumber\\
&=&\int_{\mathbb{R}^3}|u_k|^{p-2}u_k\cdot
Qu_kdx-\lambda\int_{\mathbb{R}^3}f_n(u_k)\phi_{n,u_k}\cdot
Qu_kdx+o(1)\nonumber\\
&\leq&(\int_{A_{\epsilon,k}}+\int_{\mathbb{R}^3\setminus(\varpi_k\cup
A_{\epsilon,k})}+\int_{\varpi_k})|u_k|^{p-1}\cdot
|Qu_k|dx+5\lambda
Dn^2\int_{\mathbb{R}^3}|u_k|\cdot|Qu_k|dx\nonumber\\
&\leq&\epsilon^{p-2}(\int_{A_{\epsilon,k}}|u_k|^{2}dx)^{\frac{1}{2}}(\int_{\mathbb{R}^3}|Qu_k|^{2}dx)^{1/2}+M^{p-2}(\int_{\mathbb{R}^3\setminus(\varpi_k\cup
A_{\epsilon,k})}|u_k|^{2}dx)^{\frac{1}{2}}(\int_{\mathbb{R}^3}|Qu_k|^{2}dx)^{1/2}\nonumber\\
&&+(\int_{\varpi_k}|u_k|^{p}dx)^{\frac{p-1}{p}}(\int_{\mathbb{R}^3}|Qu_k|^{p}dx)^{1/p}+5\lambda
Dn^2(\int_{\mathbb{R}^3}u^2_kdx)^{1/2}(\int_{\mathbb{R}^3}|Qu_k|^2dx)^{1/2}+o(1)\nonumber\\
&\leq&C^2\epsilon^{p-2}||u_k||^2+CM^{p-2}(\int_{\mathbb{R}^3\setminus(\varpi_k\cup
A_{\epsilon,k})}|u_k|^{2}dx)^{\frac{1}{2}}||Qu_k||
+C'(\int_{\mathbb{\varpi}_k}|u_k|^{p}dx)^{\frac{1}{p}}||Qu_k||\nonumber\\
&&+5\lambda C Dn^2||u_k||^2+o(1),
\end{eqnarray}
where $C'$ and $C$ come from (\ref{m000odiiii1}) and
(\ref{mciiiicjpp09}), respectively. Similarly, we have
\begin{eqnarray}\label{nv88ufytrtf}
||Pu_k||^2&\leq&C^2\epsilon^{p-2}||u_k||^2+CM^{p-2}(\int_{\mathbb{R}^3\setminus(\varpi_k\cup
A_{\epsilon,k})}|u_k|^{2}dx)^{\frac{1}{2}}||Pu_k||
+C'(\int_{\mathbb{\varpi}_k}|u_k|^{p}dx)^{\frac{1}{p}}||Pu_k||\nonumber\\
&&+5\lambda C Dn^2||u_k||^2+o(1).
\end{eqnarray}
Since $||u_k||^2=||Pu_k||^2+||Qu_k||^2$ (see
(\ref{nfhyfyrttr11})), these two inequalities (\ref{ijnwdvb66})
and (\ref{nv88ufytrtf}) imply that
\begin{eqnarray}\label{hfyyrtf00o}
||u_k||^2&\leq&2C^2\epsilon^{p-2}||u_k||^2+2CM^{p-2}(\int_{\mathbb{R}^3\setminus(\varpi_k\cup
A_{\epsilon,k})}|u_k|^{2}dx)^{\frac{1}{2}}||u_k||
+2C'(\int_{\mathbb{\varpi}_k}|u_k|^{p}dx)^{\frac{1}{p}}||u_k||\nonumber\\
&&+10\lambda C Dn^2||u_k||^2+o(1)
\end{eqnarray}

From $\sup_n(-\Phi_{n,\lambda}(u_k))\leq c$ and
$(1+||u_k||)||-\Phi'_{n,\lambda}(u_k)||_{X'}\rightarrow 0$, we
obtain
\begin{eqnarray}\label{nnchfdgdttdte}
o(1)+c&\geq&-\Phi_{n,\lambda}(u_{k})+\frac{1}{2}\langle
\Phi'_{n,\lambda}(u_{k}),
u_k\rangle\nonumber\\
&=&(\frac{1}{2}-\frac{1}{p})\int_{\mathbb{R}^3}|u_k|^pdx
-\frac{\lambda}{2}\int_{\mathbb{R}^3}(u_kf_n(u_k)-\frac{1}{2}F_n(u_k))\phi_{n,u_k}dx.
\end{eqnarray}
From the definitions of $f_n$ and $F_n$, we have that
\begin{eqnarray}\label{nc99dufyfff}
|t^{-2}(tf_n(t)-\frac{1}{2}F_n(t))|\leq 5/2,\ \forall
t\in\mathbb{R}.
\end{eqnarray}
From (\ref{nnchfdgdttdte}), (\ref{nc99dufyfff}), and
$\phi_{n,u_k}\leq Dn^2$, we obtain
\begin{eqnarray}\label{mmm99isgg}
&&\epsilon^{p-2}\int_{\mathbb{R}^3\setminus(\varpi_k\cup
A_{\epsilon,k})}|u_k|^2dx\nonumber\\
 &\leq&\int_{\mathbb{R}^3\setminus(\varpi_k\cup
A_{\epsilon,k})}|u_k|^pdx\nonumber\\
&\leq&c(\frac{1}{2}-\frac{1}{p})^{-1}+\frac{\lambda}{2}(\frac{1}{2}-\frac{1}{p})^{-1}\int_{\mathbb{R}^3}(u_kf_n(u_k)-\frac{1}{2}F_n(u_k))\phi_{n,u_k}dx
-\int_{\varpi_k\cup A_{\epsilon,k}}|u_k|^pdx+o(1)\nonumber\\
&\leq&c(\frac{1}{2}-\frac{1}{p})^{-1}+\frac{5\lambda}{4}(\frac{1}{2}-\frac{1}{p})^{-1}Dn^2\int_{\mathbb{R}^3}u_k^2dx+o(1)\nonumber\\
&\leq&c(\frac{1}{2}-\frac{1}{p})^{-1}+(\frac{5C^2}{4}(\frac{1}{2}-\frac{1}{p})^{-1}Dn^2)\cdot\lambda||u_k||^2+o(1).
\end{eqnarray}
It follows that
\begin{eqnarray}\label{nv99amcb66rt}
&&(\int_{\mathbb{R}^3\setminus(\varpi_k\cup
A_{\epsilon,k})}|u_k|^2dx)^{\frac{1}{2}}\nonumber\\
&\leq&
\epsilon^{1-p/2}c^{\frac{1}{2}}(\frac{1}{2}-\frac{1}{p})^{-\frac{1}{2}}
+\epsilon^{1-p/2}(\frac{5C^2}{4}(\frac{1}{2}-\frac{1}{p})^{-1}Dn^2)^{\frac{1}{2}}\lambda^{\frac{1}{2}}||u_k||+o(1).
\end{eqnarray}
If  $\epsilon$ and $\lambda'''_n$  are  such that
$$2C^2\epsilon^{p-2}=\frac{1}{8},
\
2CM^{p-2}\epsilon^{1-p/2}(\frac{5C^2}{4}(\frac{1}{2}-\frac{1}{p})^{-1}Dn^2)^{\frac{1}{2}}(\lambda'''_n)^{\frac{1}{2}}\leq\frac{1}{8},
\ 10 C Dn^2\lambda'''_n\leq\frac{1}{8},\
\lambda'''_n\leq\lambda''_n,$$ then, from (\ref{nv99amcb66rt}),
(\ref{hfyyrtf00o}) and the fact that
$\int_{\varpi_k}|u_k|^pdx=o(1)$ (Lemma \ref{hf99ofiuufyf}), we
deduce that, for $0<\lambda\leq\lambda'''_n,$ $\{||u_k||\}$ is
bounded. This completes the proof. \hfill$\Box$

\bigskip

Let $\lambda_n=\min\{\lambda'_n,\lambda''_n,\ \lambda'''_n\}$,
where $\lambda'_n$, $\lambda''_n$ and $\lambda'''_n$ are the
constants in Lemma \ref{b99d443iu}, Lemma \ref{hf99ofiuufyf} and
Lemma \ref{fbf99fifio}, respectively.

\begin{lemma}\label{nviifjgfug}
For any $n\in\mathbb{N}$ and $0<\lambda<\lambda_n$, the system
(\ref{nbgftfra44sr}) has a nontrivial solution.
\end{lemma}
\noindent{\bf Proof.} By Lemma \ref{b99d443iu}, Lemma
\ref{fbf99fifio}, and  Theorem \ref{b99d443} in the appendix, we
deduce that, for any $n\in\mathbb{N}$ and $0<\lambda<\lambda_n$,
there exists a bounded $(\overline{C})_c$ sequence $\{u_k\}$ for
$-\Phi_{n,\lambda}$ with $\inf_k|||u_k|||>0$.   Up to a
subsequence, either
\begin{itemize}
\item[{(a).}]
$\lim_{n\rightarrow\infty}\sup_{y\in\mathbb{R}^3}\int_{B_1(y)}|u_k|^2dx=0$,
or \item[{(b).}] there exist $\varrho>0$ and $a_k\in\mathbb{Z}^3$
such that $\int_{B_1(a_k)}|u_k|^2dx\geq\varrho.$
\end{itemize}
If $(a)$ occurs,  the Lions lemma (see, for example, \cite[Lemma
1.21]{Willem}) implies   $u_k$ satisfies $u_k\rightarrow 0$ in
$L^s(\mathbb{R}^3)$ for any $2<s<6$. It follows that
\begin{eqnarray}\label{nc00sowwww}
\int_{\mathbb{R}^3}|u_k|^{p-2}u_k\cdot Qu_kdx\rightarrow 0,\
\int_{\mathbb{R}^3}|u_k|^{p-2}u_k\cdot Pu_kdx\rightarrow 0,
\end{eqnarray}
and by $|t^{-1}f(t)|\leq 3$, $\forall t\in\mathbb{R}$ and Lemma
\ref{fhyr66rtfre}$(iii)$, we have
\begin{eqnarray}\label{n00dfdhfgfff}
\int_{\mathbb{R}^3}f_n(u_k)\phi_{n,u_k}\cdot Qu_kdx
&=&\int_{\mathbb{R}^3}u^{-1}_kf_n(u_k)\phi_{n,u_k}\cdot u_k\cdot
Qu_kdx\nonumber\\
&\leq&3\int_{\mathbb{R}^3}u_k\phi_{n,u_k} \cdot
Qu_kdx\nonumber\\
&\leq&3(\int_{\mathbb{R}^3}|\phi_{n, u_k}|^{6}dx)^{\frac{1}{6}}
(\int_{\mathbb{R}^3}|u_k|^{\frac{12}{5}}dx)^{\frac{5}{6}}(\int_{\mathbb{R}^3}|Qu_k|^{\frac{12}{5}}dx)^{\frac{5}{6}}\nonumber\\
&\leq&
3C^{-1}_*(\int_{\mathbb{R}^3}|u_k|^{\frac{12}{5}}dx)^{\frac{5}{3}}(\int_{\mathbb{R}^3}|Qu_k|^{\frac{12}{5}}dx)^{\frac{5}{6}}\rightarrow
0.
\end{eqnarray}
Similarly, we have $\int_{\mathbb{R}^3}f_n(u_k)\phi_{n,u_k}\cdot
Pu_kdx\rightarrow 0.$ Then from (\ref{nc00sowwww}),
(\ref{n00dfdhfgfff}), (\ref{hf999f876r6rr}) and
(\ref{hfuu7g7tyyt1}), we obtain $||u_k||\rightarrow 0.$ This
contradicts $\inf_{k}|||u_k|||>0$. Therefore, case $(a)$ cannot
occur. As case  $(b)$ therefore occurs, $w_k=u_k(\cdot+a_k)$
satisfies $w_k\rightharpoonup u_0\neq 0$. From
$(1+||w_k||)||-\Phi'_{n,\lambda}(w_k)||_{X'}=(1+||u_k||)||-\Phi'_{n,\lambda}(u_k)||_{X'}\rightarrow
0$ and the weakly sequential continuity of $-\Phi'_{n,\lambda}$
(see Lemma \ref{b99d443iu}), we have that
$-\Phi'_{n,\lambda}(u_0)=0.$ Therefore, $(u_0,\phi_{n,u_0})$ is a
nontrivial solution of
 (\ref{nbgftfra44sr}). This completes the proof.\hfill$\Box$

\section{Proof of Theorem \ref{th1}}

\begin{proposition}\label{nc99dyftrf}
Suppose that $3<p<6$. Then for any $\Lambda>0$, there exists
$N_\Lambda>0$ such that, if $(u,\phi)\in H^1(\mathbb{R}^3)\times
\mathcal{D}^{1,2}(\mathbb{R}^3)$ is a solution of
(\ref{nbgftfra44sr}) with  $n\geq N_\Lambda$ and
$0<\lambda\leq\Lambda$, then
$$|u(x)|\leq n,\ \forall x\in\mathbb{R}^3.$$
\end{proposition}
\noindent{\bf Proof.} We apply an indirect argument, and assume by
contradiction that there exist $\Lambda_0>0$, a real number
sequence $\{\mu_n\}$  and a sequence $\{(u_n,\phi_n)\}$ in $
H^1(\mathbb{R}^3)\times \mathcal{D}^{1,2}(\mathbb{R}^3)$ such that
$0<\mu_n\leq \Lambda_0$,  $(u_n,\phi_n)$ is a solution of
(\ref{nbgftfra44sr}) with $\lambda=\mu_n$  and
\begin{eqnarray}\label{nciicjhfii}
||u_n||_{L^\infty(\mathbb{R}^3)}>n.
\end{eqnarray}

Since $\phi_{n,u}$ is a bounded function in $\mathbb{R}^3$ (see
Lemma \ref{fhyr66rtfre}$(iv)$)  and $u_n\in H^1(\mathbb{R}^3)$ is
a solution of
\begin{eqnarray}\label{bcnhhfyftff}
-\Delta u+V(x)u+|u|^{p-2}u=\mu_nf_n(u)\phi_{n,u}\quad \mbox{in}\
\mathbb{R}^3,
\end{eqnarray}
 the bootstrap argument of  elliptic equations (see \cite{GT}),
implies that $u_n\in C^1(\mathbb{R}^3)$. Let $x_n\in\mathbb{R}^3$
be such that
 \begin{eqnarray}\label{bv77tffvb}
|
u_n(x_n)|=\max_{x\in\mathbb{R}^3}|u_n(x)|=||u_n||_{L^\infty(\mathbb{R}^3)}.
 \end{eqnarray}

We shall use the blow-up argument of \cite{bahri} to induce a
contradiction. Let $M_n=\max_{x\in\mathbb{R}^3}|u_n(x)|$ and
$$\tilde{u}_n(x)=M^{-1}_n u_n(x_n+xM^{-\alpha}_{n}),\ x\in\mathbb{R}^3$$
with $\alpha=(p-2)/2$. Then,
\begin{eqnarray}\label{nvmiivjfh}
|\tilde{u}_n|\leq 1\ \mbox{in}\ \mathbb{R}^3,\
|\tilde{u}_n(0)|=1,\ \forall n\in\mathbb{N},
\end{eqnarray} and
$$u_n(x)=M_n\tilde{u}_n(M^\alpha_n(x-x_n)),\ x\in\mathbb{R}^3.$$
Substituting this expression into (\ref{bcnhhfyftff}), we obtain
\begin{eqnarray}\label{nc99cifjfh}
&&-M^{1+2\alpha}_n(\Delta
\tilde{u}_n)(M^\alpha_n(x-x_n))+M_nV(x)\tilde{u}_n(M^\alpha_n(x-x_n))\nonumber\\
&&+M^{p-1}_n(|\tilde{u}_n|^{p-2}\tilde{u}_n)(M^\alpha_n(x-x_n))\nonumber\\
&=&\mu_nf_n(M_n\tilde{u}_n(M^\alpha_n(x-x_n)))\int_{\mathbb{R}^3}\frac{F_n(M_n\tilde{u}_n(M^\alpha_n(y-x_n)))}{4\pi|x-y|}dy.
\end{eqnarray}
By a direct computation, we have
\begin{eqnarray}\label{mvc00voigu}
\int_{\mathbb{R}^3}\frac{F_n(M_n\tilde{u}_n(M^\alpha_n(y-x_n)))}{4\pi|x-y|}dy=
M^{-2\alpha}_n\int_{\mathbb{R}^3}\frac{F_n(M_n\tilde{u}_n(y))}{4\pi|M^\alpha_n(x-x_n)-y|}dy,\nonumber
\end{eqnarray}
Together with (\ref{nc99cifjfh}), this yields
\begin{eqnarray}\label{mc00ovifyu}
&&-\Delta
\tilde{u}_n+M^{2-p}_nV(x_n+xM^{-\alpha}_n)\tilde{u}_n+|\tilde{u}_n|^{p-2}\tilde{u}_n\nonumber\\
&=&M^{-(2p-3)}_n\cdot\mu_nf_n(M_n\tilde{u}_n)
\int_{\mathbb{R}^3}\frac{F_n(M_n\tilde{u}_n(y))}{4\pi|x-y|}dy.
\end{eqnarray}
From the definition of $f_n$, we have $|f_n(t)|\leq 6n$, $\forall
t\in\mathbb{R}$. Then,   Lemma \ref{fhyr66rtfre}$(iv)$ gives
$$0\leq\phi_{n,M_n\tilde{u}_n}=\int_{\mathbb{R}^3}\frac{F_n(M_n\tilde{u}_n(y))}{4\pi|x-y|}dy\leq Dn^2.$$
Then by $p>3$,  $\mu_n\leq \Lambda_0$ and  $M_n>n$, we obtain
\begin{eqnarray}\label{mnuuyhtgq}
M^{-(2p-3)}_n\cdot\mu_n|f_n(M_n\tilde{u}_n)|
\int_{\mathbb{R}^3}\frac{F_n(M_n\tilde{u}_n(y))}{4\pi|x-y|}dy\leq
6D\Lambda_0n^3M^{-(2p-3)}_n\rightarrow0,\ n\rightarrow\infty.
\end{eqnarray}
Moreover, because $V$ is a bounded function in $\mathbb{R}^3$, we
have that
\begin{eqnarray}\label{mvoovigfufj}
M^{2-p}_nV(x_n+xM^{-\alpha}_n)\rightarrow0,\ n\rightarrow\infty
\end{eqnarray} holds uniformly for
$x\in\mathbb{R}^3$. Then, by the standard elliptic estimates (see
Section 9.2 of \cite{GT}), we deduce from $|\tilde{u}_n|\leq 1$ in
$ \mathbb{R}^3,$ (\ref{mc00ovifyu}), (\ref{mnuuyhtgq}) and
(\ref{mvoovigfufj}) that, for any $2\leq q<\infty$ and
$0<R<\infty,$ $\tilde{u}_n$ is bounded in $W^{2,q}(B_R(0))$.

Without loss of generality, we may assume that $\tilde{u}_n$
converges weakly in $W^{2,q} (B_R(0))$ ($\forall R <+\infty$, $
\forall p <+\infty$) and thus in particular in $C^1 (B_R(0))$ to
some $u_0$ satisfying $|u_0(0)|=1$. From  (\ref{mc00ovifyu}),
(\ref{mnuuyhtgq})
 and (\ref{mvoovigfufj}), we can see that $u_0$
satisfies
\begin{eqnarray}\label{hcgttdrd}
-\Delta u+|u|^{p-2}u=0,\ u\in C^1(\mathbb{R}^3).
\end{eqnarray}
However,  by Theorem 1 of \cite{brezis}, the only solution to this
equation is $u=0,$ which  contradicts $|u_0(0)|=1$. This completes
the proof.\hfill$\Box$

\bigskip

\noindent{\bf Proof of Theorem \ref{th1}.} We choose $\Lambda=1$
in Proposition \ref{nc99dyftrf} and choose $n_*\in \mathbb{N}$
satisfying $n_*\geq N_\Lambda$. Then by Lemma \ref{nviifjgfug} and
Proposition \ref{nc99dyftrf},  for any
$0<\lambda<\lambda_0:=\min\{\lambda_{n_*},1\}$,  problem
(\ref{nbgftfra44sr}) with $n=n_*$ has a nontrivial  solution
$(u,\phi)$ satisfying
$$||u||_{L^\infty(\mathbb{R}^3)}\leq n_*.$$
It follows that $f_{n_*}(u)=u$ and $F_{n_*}(u)=u^2.$ Hence,
$(u,\phi)$ is a nontrivial solution of
(\ref{nbgftfrr}).\hfill$\Box$

\section{Appendix: A variant infinite-dimensional linking theorem}\label{nvb6ftrf}

In this section, we give a new infinite-dimensional linking
theorem.  This theorem replaces the $\tau$-upper semi-continuous
assumption (see (6.3) in \cite{Willem}) in the Kryszewski and
Szulkin's infinite-dimensional linking theorem (see \cite[Theorem
6.10]{Willem} or \cite[Theorem 3.4]{KS}) with other assumptions
(see (\ref{nx99s8s7yy}) in Theorem \ref{b99d443}). Our theorem is
a generalization of the classical finite-dimensional linking
theorem (see \cite[Theorem 5.3]{rabinowitz})

Before state the infinite-dimensional linking  theorem, we give
some notations and definitions.

Let $X$ be a separable Hilbert space with inner product
$(\cdot,\cdot)$ and norm $||\cdot||$, respectively. $Y$ and $Z$
are closed subspaces of $X$ and $X=Y\oplus Z.$ Let $\{e_k\}$ be a
total orthonormal sequence in $Y$. Let
\begin{eqnarray}\label{nchu777d6d66}
Q:X\rightarrow Z,\ P:X\rightarrow Y
\end{eqnarray} be the
orthogonal projections. We define
\begin{eqnarray}\label{o88dytgdg}
|||u|||=\max\Big\{||Qu||,\sum^{\infty}_{k=1}\frac{1}{2^{k+1}}|(Pu,e_k)|\Big\}
\end{eqnarray}
on $X.$ Then, $$||Qu||\leq|||u|||\leq ||u||,\ \forall u\in X,$$
and  if $||u_n||$ is bounded and $|||u_n-u|||\rightarrow 0$, then
$\{u_n\}$ weakly converges to $u$ in $X.$
 The topology
generated by $|||\cdot|||$ is denoted by $\tau$, and all
topological notations related to it will include this symbol.

\medskip

Let $R > r > 0$ and  $u_0 \in Z$ with $||u_0||
 = 1$.
 Set
\begin{eqnarray}\label{ncg7777qwwq}
N = \{u \in Z \ |\  ||u||
 = r\},\  M = \{u+tu_0\ |\  u\in Y ,\ t\geq 0,\
||u+tu_0|| \leq R\}. \end{eqnarray} Then, $M$ is a submanifold of
$Y \oplus \mathbb{R}^+u_0$ with boundary
\begin{eqnarray}\label{m9155sfddd}
\partial M=\{u\in Y\ |\ ||u||\leq
R\}\cup\{u+tu_0\ |\ u\in Y,\ t>0,\ ||u+tu_0||=R \}.
\end{eqnarray}

\begin{definition}\label{bcvujnbytg}
Let  $J\in  C^1(X, \mathbb{R})$. A sequence $\{u_n\} \subset X$ is
called a  $(\overline{C})_c$ sequence  for $J$, if $$ \sup_nJ(u_n)
\leq c\quad \mbox{ and}\quad
(1+||u_n||)||J'(u_n)||_{X'}\rightarrow 0,\ \mbox{ as}\
n\rightarrow\infty.$$
\end{definition}

\begin{Theorem} \label{b99d443}  If
 $H\in C^1(X,\mathbb{R})$ satisfies
\begin{description}  \item {$(a)$}  $H'$  is weakly
sequentially continuous, i.e., if $u\in X$ and $\{u_n\}\subset X$
are such that $u_n\rightharpoonup u$, then, for any $\varphi\in
X,$ $\langle H'(u_n),\varphi\rangle\rightarrow \langle
H'(u),\varphi\rangle$.  \item {$(b)$} there exist $\delta>0,$ $u_0
\in Z  $ with $||u_0||=1$, and $R
> r
> 0$ such that
\begin{eqnarray}\label{nx99s8s7yy}
\inf_N H>\max\Big\{\sup_{\partial M} H,
\sup_{|||u|||\leq\delta}H(u)\Big\}
\end{eqnarray}
and
\begin{eqnarray}\label{m99iuxyyxy}
\sup_MH<+\infty,
\end{eqnarray}
\end{description}  Then there exists a $(\overline{C})_c$ sequence $\{u_n\}$ for $H$ with $c=\sup_MH$ and
 $\inf_{n}|||u_n|||\geq\delta/2$.
\end{Theorem}
\noindent{\bf Proof.} Arguing indirectly, assume that  the result
does not hold. Then, there exists $\epsilon>0$ such that
\begin{eqnarray}\label{nchyyftdg}
(1+||u||)||H'(u)||_{X'}\geq\epsilon,\ \forall u\in E
\end{eqnarray}
where $$E=\{u\in X\ |\ H(u)\leq d+\epsilon\}\cap\{u\in X\ |\
|||u|||\geq\delta/2\}$$ and $$d=\sup_{M}H.$$ From
(\ref{nx99s8s7yy}), we can choose $\epsilon$ such that
\begin{eqnarray}\label{nc88c7d6ttd}
0<\epsilon<\inf_NH-\max\Big\{\sup_{\partial M} H,
\sup_{|||u|||\leq\delta}H(u)\Big\}.
\end{eqnarray}

$\bf{Step 1.}$ A  vector field in a $\tau$-neighborhood of $E$.

Let
\begin{eqnarray}\label{hvcnbgdtdt66t11}
b=\inf_N H,\quad T=2(d-b+2\epsilon)/\epsilon,\quad R=(1+\sup_{u\in
M}||u||)e^T
\end{eqnarray}
and
\begin{eqnarray}\label{hhbc777cyttdg}
B_R=\{u\in X\ |\ ||u||\leq R\}.
\end{eqnarray}

For every $u\in E\cap B_R,$  there exists $\phi_u\in X$ with
$||\phi_u||=1$ such that $\langle
H'(u),\phi_u\rangle\geq\frac{3}{4}||H'(u)||_{X'}$. Then,
(\ref{nchyyftdg}) implies
\begin{eqnarray}\label{n99v8f77fy}
(1+||u||)\langle H'(u),\phi_u\rangle>\frac{1}{2}\epsilon.
\end{eqnarray}
From the definition of $|||\cdot|||$, we deduce that  if a
sequence $\{u_n\}\subset  E\cap B_R$  $\tau$-converges  to $u\in
X$, i.e., $|||u_n-u|||\rightarrow 0$, then $u_n\rightharpoonup u$
in $X$ (see Remark 6.1 of \cite{Willem}).   By the weakly
sequential continuity of $H'$, we get that for any $\varphi\in X$,
$\langle H'(u_n),\varphi\rangle\rightarrow\langle
H'(u),\varphi\rangle$. This implies that $H'$ is
$\tau$-sequentially continuous in $E\cap B_R$. By
(\ref{n99v8f77fy}), the $\tau$-sequential continuity of $H'$ in
$E\cap B_R$ and the weakly lower semi-continuity of the norm
$||\cdot||$, we get that  there exists a $\tau$-open neighborhood
$V_u$ of $u$ such that
\begin{eqnarray}\label{nv88vuijf}
\langle H'(v),(1+||u||)\phi_u\rangle>\frac{1}{2}\epsilon,\ \forall
v\in V_u,
\end{eqnarray}
and
\begin{eqnarray}\label{nncbyyft55er}
||(1+||u||)\phi_u||=1+||u||\leq 2(1+||v||),\ \forall v\in V_u.
\end{eqnarray}

Because $B_R$ is a bounded convex closed set in the Hilbert space
$X$, $B_R$ is a $\tau$-closed set. Therefore, $X\setminus B_R$ is
a $\tau$-open set.

 The family
$$\mathcal{N}=\{V_u\ |\ u\in E\cap B_R\}\cup\{X\setminus B_R\}$$
is a $\tau$-open covering of $E$. Let
$$\mathcal{V}=\Big(\bigcup_{u\in E\cap B_R} V_u\Big)\bigcup (X\setminus B_R).$$
Then, $\mathcal{V}$ is a $\tau$-open neighborhood of $E.$

Since $\mathcal{V}$ is metric, hence paracompact, there exists a
local finite $\tau$-open covering $\mathcal{M}=\{M_i\ |\ i\in
\Lambda\}$ of $\mathcal{V}$ finer than $\mathcal{N}$. If
$M_i\subset V_{u_i}$ for some $u_i\in E$, we choose
$\varpi_i=(1+||u_i||)\phi_{u_i}$ and if $M_i\subset X\setminus
B_R$, we choose $\varpi_i=0$.
 Let
$\{\lambda_i\ |\ i\in I\}$ be a $\tau$-Lipschitz continuous
partition of unity subordinated to $\mathcal{M}$. And let
\begin{eqnarray}\label{ncb77dtdf}
\xi(u):=\sum_{i\in I}\lambda_i(u)\varpi_i,\ u\in \mathcal{V}.
\end{eqnarray}
Since the  $\tau$-open covering $\mathcal{M}$ of $\mathcal{V}$ is
local finite,  each $u\in\mathcal{V}$ belongs to only finite many
sets $M_i$. Therefore,  for every  $u\in\mathcal{V}$, the sum in
(\ref{ncb77dtdf}) is only a finite sum. It follows that, for any
$u\in \mathcal{V},$ there exist a $\tau$-open neighborhood
$U_u\subset\mathcal{V}$ of $u$ and $L_u>0$ such that
 $\xi(U_u)$ is
contained in a finite-dimensional subspace of $X$ and
\begin{eqnarray}\label{mm99ifuyy}
||\xi(v)-\xi(w)||\leq L_u|||v-w|||,\ \forall v,w\in U_u.
\end{eqnarray}
Moreover,  by the definition of $\xi$, (\ref{nv88vuijf}) and
(\ref{nncbyyft55er}), we get that, for every $u\in \mathcal{V},$
\begin{eqnarray}\label{oikjnhh7yyt}
||\xi(u)||\leq 1+||u||\quad \mbox{and}\quad \langle H'(u),
\xi(u)\rangle\geq 0
\end{eqnarray}
and for every $u\in E\cap B_R,$
\begin{eqnarray}\label{hhcbtgdtdr}
\langle H'(u), \xi(u)\rangle> \frac{1}{2}\epsilon.
\end{eqnarray}

$\bf{Step 2.} $  Let $\theta$ be a smooth  function satisfying
$0\leq \theta\leq 1$ in $\mathbb{R}$ and
$$\theta(t)=\left\{
\begin{array}
[c]{ll}
0 , & t\leq \frac{2\delta}{3},\\
1, & t\geq \delta.
\end{array}
\right.\label{e1}$$ Let
$$
\chi(u)=\left\{\begin{array} [c]{ll}
-\theta(|||u|||)\xi(u),&\ u\in\mathcal{V}\ ,  \\
0, & |||u||| \leq \frac{2\delta}{3}. \end{array}\right. \label{e1}
$$ Then, $\chi$ is a vector field defined in
$$\mathcal{W}=\mathcal{V}\cup\{u\in X\ |\ |||u|||<\delta\}.$$
It is a $\tau$-open neighborhood of $H^{d+\epsilon}\cup(X\setminus
B_R)$, where
$$H^{d+\epsilon}:=\{u\in X\ |\ H(u)\leq d+\epsilon\}.$$ From
(\ref{mm99ifuyy}), (\ref{hhcbtgdtdr}) and the definition of
$\chi$,  we deduce that the mapping $\chi$ satisfies that
\begin{description}
\item{$(\bf{a}).$} each $u\in \mathcal{W}$ has a $\tau$-open set
$V_u$ such that $\chi(V_u)$ is contained in a finite-dimensional
subspace of $X,$

\item{$(\bf{b}).$} for any $u\in \mathcal{W},$ there exist a
$\tau$-open neighborhood  $U_u$ of $u$ and $L'_u>0$ such that
\begin{eqnarray}\label{m000oiu76}
||\chi(v)-\chi(w)||\leq L'_u|||v-w|||,\ \forall v,w\in U_u.
\end{eqnarray}
This means that $\chi$ is locally Lipschitz continuous and
$\tau$-locally Lipschitz continuous,

\item{$(\bf{c}).$}
\begin{eqnarray}\label{bcvyyftfg} ||\chi(u)||\leq 1+||u||,\ \forall u\in
\mathcal{W},\end{eqnarray}

\item{$(\bf{d}).$}
\begin{eqnarray}\label{b7dyttdra} \langle
H'(u),\chi(u)\rangle\leq 0,\ \forall u\in \mathcal{W}.
\end{eqnarray} and
\begin{eqnarray}\label{b7dyttdratr} \langle
H'(u),\chi(u)\rangle<-\frac{1}{2}\epsilon,\ \forall u\in \{u\in E\
|\ |||u|||\geq\delta\}\cap B_R.
\end{eqnarray}
\end{description}

$\bf{Step 3.} $ From (\ref{m000oiu76}) and the fact that
$|||v|||\leq ||v||$, $\forall v\in X$, we have
$$
||\chi(v)-\chi(w)||\leq L'_u||v-w||,\ \forall v,w\in U_u. $$ This
implies that $\chi$ is a local Lipschitz mapping under the
$||\cdot||$ norm. Then by the standard theory of ordinary
differential equation in Banach space,  we deduce that the
following initial value problem
  \begin{eqnarray}\label{iijuytytyyy16}
\left\{
\begin{array}
[c]{ll}
\frac{d\eta}{dt}=\chi(\eta),\\
\eta(0,u)=u\in \mathcal{W}.
\end{array}
\right.
\end{eqnarray}
has a unique solution in $\mathcal{W}$ , denoted by $\eta(t, u)$,
with right maximal interval of existence $[0, T(u))$.
 Furthermore, using
(\ref{m000oiu76}) and the Gronwall inequality (see, for example,
Lemma 6.9 of \cite{Willem}), the similar argument  as the proof of
$\bf c)$ in \cite[Lemma 6.8]{Willem}  yields that
\begin{description}
\item{$(\bf{A}).$} $\eta$ is $\tau$-continuous, i.e., if
 $u_n\in \mathcal{W}$,  $u_0\in
\mathcal{W}$, $0\leq t_n<T(u_n)$ and $0\leq t_0<T(u_0)$ satisfy
$|||u_n-u_0|||\rightarrow 0$ and $t_n\rightarrow t_0$, then
$|||\eta(t_n,u_n)-\eta(t_0,u_0)|||\rightarrow 0.$
\end{description}

From (\ref{b7dyttdra}), we have $$\frac{d}{dt}H(\eta(t,u))=\langle
H'(\eta(t,u)),\eta(t,u)\rangle\leq0.$$  Therefore, $H$ is
non-increasing along the flow $\eta.$  It follows that
$\{\eta(t,u)\ |\ 0\leq t\leq T(u)\}\subset H^{d+\epsilon}$ if
$u\in H^{d+\epsilon},$ i.e., $H^{d+\epsilon}$ is an invariant set
of the flow $\eta.$ Then, (\ref{bcvyyftfg}) and Theorem 5.6.1 of
\cite{laksh} (or Corollary 4.6 of \cite{schechter}) implies  that,
for any $ u\in H^{d+\epsilon},$
   $T(u)=+\infty$.

   \medskip

{\bf Step 4. } We shall prove that
\begin{eqnarray}\label{mxjvchhvc99987}
\{\eta(t, u)\ |\ 0\leq t\leq T,\ u\in M\}\subset B_R.
\end{eqnarray}

Let $u\in H^{d+\epsilon}.$ By the result in Step 3,  we have
$T(u)=+\infty$ and
$$\eta(t,u)=u+\int^t_0 \chi(\eta(s,u))ds,\quad \forall t\in[0,+\infty).$$
Together with (\ref{bcvyyftfg}), this yields
\begin{eqnarray}\label{ufghyyft88u6}
||\eta(t,u)||\leq ||u||+\int^t_0 ||\chi(\eta(s,u))||ds\leq
||u||+\int^t_0(1+ ||\eta(s,u))||)ds.\nonumber
\end{eqnarray}
Then, by the Gronwall inequality (see, for example, Lemma 6.9 of
\cite{Willem}), we get that
\begin{eqnarray}\label{jfuudt5555}
||\eta(t,u)||\leq (1+||u||)e^t-1, \quad \forall t\in[0,+\infty).
\end{eqnarray}
Since $M\subset H^{d+\epsilon}$, by (\ref{jfuudt5555}) and the
definition of $R$ (see \ref{hvcnbgdtdt66t11}), we get
(\ref{mxjvchhvc99987}).

\medskip

{\bf Step 5. } From  the choice of $\epsilon$ (see
(\ref{nc88c7d6ttd})), we have
$$\sup_{|||u|||\leq\delta}H<b-\epsilon.$$
It follows that
$$\{u\in X\ |\ |||u|||\leq\delta\}\subset H^{b-\epsilon}
:=\{u\in X\ |\ H(u)\leq b-\epsilon\}. $$ Together with
(\ref{b7dyttdratr}), this  yields
\begin{eqnarray}\label{hhdyttd}
\langle H'(u),\chi(u)\rangle<-\frac{1}{2}\epsilon,\ \forall u\in
H^{d+\epsilon}_{b-\epsilon}\cap B_R, \end{eqnarray} where
$$H^{d+\epsilon}_{b-\epsilon}:= \{u\in X\ |\ b-\epsilon\leq
H(u)\leq d+\epsilon\}.$$

We show that, for any $u\in M$, $H(\eta(T,u))\leq b-\epsilon.$
Arguing indirectly, assume that this were not true.  Then,  there
exists $u\in M$ such that $H(\eta(T,u))> b-\epsilon.$ Since $H$ is
non-increasing along the flow $\eta$, from (\ref{mxjvchhvc99987}),
we deduce that  $\{\eta(t,u)\ |\ 0\leq t\leq T\}\subset
H^{d+\epsilon}_{b-\epsilon}\cap B_R$. Then, by (\ref{hhdyttd}),
\begin{eqnarray}\label{nc88cudyd}
H(\eta(T,u))&=&H(\eta(0,u))+\int^T_0\Big\langle
H'(\eta(s,u)),\chi(\eta(s,u))\Big
\rangle ds\nonumber\\
&\leq&H(\eta(0,u))+\int^T_0(-\frac{1}{2}\epsilon)ds\nonumber\\
&\leq&d+\epsilon-\frac{1}{2}\epsilon T=b-\epsilon.
\end{eqnarray}
This contradicts $H(\eta(T,u))>b-\epsilon$. Therefore, we have
 \begin{description}
\item{$(\bf{B}).$} $\eta(T, M)\subset H^{b-\epsilon}.$
\end{description}
Moreover, using the result $\bf (a)$ in Step 2 and  and the fact
that $\eta$ is $\tau$-continuous (see $\bf (A)$), the similar
argument as the proof of the result $b)$  of \cite[Lemma
6.8]{Willem} yields that
\begin{description}
\item{$(\bf{C}).$} Each point $(t,u)\in [0,T]\times
H^{d+\epsilon}$ has a $\tau$-neighborhood $N_{(t,u)}$ such that
$$\{v-\eta(s,v)\ |\ (s,v)\in N_{(t,u)}\cap ([0,T]\times
H^{d+\epsilon})\}$$ is contained in a finite-dimensional subspace
of $X.$
\end{description}

$\bf{Step 6.}$  Let $$h:[0,T]\times M\rightarrow X,\
h(t,u)=P\eta(t,u)+(||Q\eta(t,u)||-r)u_0$$ where $P, Q$, $r$ and
$u_0$ are defined in (\ref{nchu777d6d66}), (\ref{ncg7777qwwq}) and
(\ref{m9155sfddd}). Then $$0\in h(t,M)\Leftrightarrow
\eta(t,M)\cap N\neq\emptyset.$$ From $\inf_NH>\sup_{\partial M} H$
(see (\ref{nx99s8s7yy})) and the fact that, for any $u\in X$, the
function $H(\eta(\cdot, u))$ is non-increasing, we deduce that
$\inf_NH>\sup_{u\in\partial M}H(\eta(t,u))$, $\forall t\in [0,T]$.
Therefore,
\begin{eqnarray}\label{tgyyyuujpoik}0\not\in h(t,\partial M), \ \forall t\in
[0,T]. \end{eqnarray}  Since $\eta$ has the properties $\bf (A)$
and $\bf (C)$ obtained in step 3 and step 5 respectively and $h$
satisfies (\ref{tgyyyuujpoik}), there is an appropriate degree
theory for $\deg(h(t,\cdot), M,0)$ (see Proposition 6.4 and
Theorem 6.6 of \cite{Willem}).  Then, the same argument as the
proof of Theorem 6.10 of \cite{Willem} yields that
$$\deg(h(T,\cdot), M, 0)=\deg(h(0,\cdot), M, 0)\neq 0.$$
It follows that $0\in h(T,M)$ and $\eta(T,M)\cap N\neq\emptyset.$
Therefore, there exists $u\in M$ such that $H(\eta(T,u))\geq b.$
It contradicts the property $\bf (B)$ obtained  in step 5. This
completes the proof of this theorem.\hfill$\Box$

\end{document}